\documentclass[a4paper,12pt,final]{article}
\usepackage{amsfonts,amsmath}
\usepackage[margin=1in]{geometry}
\usepackage[colorlinks=true,linkcolor=black,citecolor=blue,urlcolor=blue,]{hyperref}
\usepackage[numbers]{natbib}
\usepackage{mathtools}
\usepackage{booktabs}

\title{\Large Dual-Domain Deep Learning Method to Accelerate Local Basis Functions Computation for Reservoir Simulation in High-Contrast Heterogeneous Porous Media\footnote{This study was supported by the Postgraduate Research Scholarship of Xi'an Jiaotong-Liverpool University (Grant No. FOSA2412003).}}

\author{Peiqi Li\thanks{School of Mathematics and Physics, Xi'an Jiaotong-Liverpool University, 215123, Suzhou, Jiangsu, China (Email: LPQ\_0619@outlook.com).} 
\and
Jie Chen*\thanks{{School of Mathematics and Physics, Xi'an Jiaotong-Liverpool University, 215123, Suzhou, Jiangsu, China (Email: Jie.Chen01@xjtlu.edu.cn).}}}

\begin{document}
\date{}
\maketitle

\begin{abstract}
In energy science, Darcy flow in heterogeneous porous media is a central problem in reservoir simulation. However, the pronounced multiscale characteristics of such media pose significant challenges to conventional numerical methods in terms of computational demand and efficiency. The Mixed Generalized Multiscale Finite Element Method (MGMsFEM) provides an effective framework for addressing these challenges, yet the construction of multiscale basis functions remains computationally expensive. In this work, we propose a dual-domain deep learning framework to accelerate the computation of multiscale basis functions within MGMsFEM for solving Darcy flow problems. By extracting and decoding permeability field features in both the frequency and spatial domains, the method enables rapid generation of numerical matrices of multiscale basis functions. Numerical experiments demonstrate that the proposed framework achieves significant computational acceleration while maintaining high approximation accuracy, thereby offering the potential for future applications in real-world reservoir engineering.
\end{abstract}


\section{Introduction}
\label{sec:intro}

Reservoir simulation occupies a central position in energy science and engineering. It not only provides theoretical and technical support for the exploration and development of oil and gas resources, but also plays an important role in geological carbon sequestration and groundwater resource man-agement\cite{chen2006computational,lake1989enhanced}. With the ongoing energy transition and the increasing severity of environmental issues, high-accuracy and high-efficiency reservoir simulation methods are of great significance for enhancing resource recovery and ensuring energy security.

In reservoir simulation, the Darcy flow model serves as the fundamental mathematical framework for describing fluid motion in porous media and constitutes the theoretical basis for constructing multiphase and multicomponent flow models\cite{durlofsky1991numerical}. However, subsurface porous media often exhibit pronounced heterogeneity, with permeability distributions spanning several orders of magnitude, leading to complex nonuniformity and strong multiscale characteristics.

These multiscale characteristics pose significant challenges for numerical computation. On the one hand, employing fine-grid discretization to solve the Darcy flow problem can effectively capture fine-scale features, but it leads to a dramatic increase in computational scale and extremely high computational cost. On the other hand, coarse-grid approximations can reduce the computational burden, yet they often sacrifice critical fine-scale information, resulting in a loss of accuracy\cite{jenny2003multi}. Therefore, achieving high computational efficiency while maintaining accuracy has become a critical issue in reservoir simulation.

To address the multiscale challenges of Darcy flow in heterogeneous porous media, various nu-merical methods have been developed. Traditional approaches such as the finite element method\cite{wheeler2006multipoint} (FEM), the finite volume method\cite{eymard2000finite} (FVM), and the finite difference method (FDM) perform well in homogeneous or mildly heterogeneous media, but they often become inefficient when applied to strongly heterogeneous and multiscale problems\cite{jenny2003multi}. 

In recent years, multiscale numerical methods have emerged as powerful tools for tackling such challenges. These methods construct multiscale basis functions on coarse grids that encode fine-scale information, thereby significantly reducing computational cost while maintaining high accuracy. Representative approaches include the multiscale volume method\cite{hajibeygi2008iterative,jenny2005adaptive} (MsFVM), the multiscale finite element method\cite{chen2003mixed} (MsFEM), the generalized multiscale finite element method\cite{efendiev2013generalized} (GMsFEM), and the mixed generalized multiscale finite element method\cite{chung2015mixed,chen2020generalized} (MGMsFEM). Among them, the MGMsFEM has received particular attention in reservoir simulation and subsurface flow problems due to its ability to preserve local mass conservation while ensuring global accuracy.

However, existing multiscale methods still suffer from significant limitations. Taking the MGMsFEM as an example, its core computational step lies in the construction of multiscale basis functions, which typically requires solving local eigenvalue problems or boundary problems in each coarse block. This process is computationally expensive, especially in large-scale three-dimensional reservoir simulations (such as SPE10 benchmark model\cite{christie2001tenth}), where enormous number of basis functions and repeated solutions of local problems severely restrict overall efficiency. Therefore, achieving substantial acceleration of basis function construction while preserving the high-accuracy characteristics of MGMsFEM has become a critical challenge for advancing its applicability in real-world engineering practice.

The rapid development of deep learning has opened new avenues for accelerating scientific compu-ting in multiscale porous media problems\cite{choubineh2022innovative,liu2024mitigating,choubineh2023deep,geng2024swin,liu2024learning}. The Fourier Neural Operator (FNO) is capable of learning mappings between function spaces and has demonstrated remarkable efficiency in solving partial differential equations\cite{li2020fourier,li2023fourier} (PDEs). Physics-Informed Neural Networks (PINNs) embed governing equations directly into the training process, thereby enforcing physical constraints and enabling solution approximation without requiring extensive data\cite{raissi2019physics}. However, this approach shows clear limitations when dealing with more complex equations and in accurately fitting boundary conditions. Similarly, Deep Operator Networks (DeepONets) provide a flexible framework for ap-proximating nonlinear operators and exhibit strong generalization capabilities\cite{lu2021learning}. More recently, Kolmogorov–Arnold Networks (KANs) have been proposed as a novel and interpretable architec-ture, showing promise in tackling high-dimensional problems\cite{liu2024kan}. Despite their differences, these architectures share the common advantage of effectively capturing multiscale features while re-ducing computational cost.

Nevertheless, their applications in reservoir simulation have primarily focused on approximating global solutions, while relatively little attention has been devoted to accelerating the computation of multiscale basis functions within frameworks such as the MGMsFEM. Meanwhile, in the field of energy science, reservoir simulation plays a critical role in enhancing the efficient recovery of oil and gas, enabling geological carbon sequestration, and supporting optimized groundwater resource management, all of which place increasingly stringent demands on the accuracy and efficiency of numerical methods. If deep learning can be organically integrated with multiscale numerical methods to achieve substantial acceleration while maintaining high approximation accuracy, it will not only help overcome the computational bottlenecks of current numerical simulations but also provide practical and feasible pathways for real-world engineering applications in energy science.

In this study, we attempt to integrate multiscale numerical methods with deep learning to accelerate the computation of multiscale basis functions in MGMsFEM. Building upon the FNO and convolutional kernels of varying sizes, we perform feature extraction in both the frequency and spatial domains. The extracted features from different kernel sizes are then fused through an additive operator to enable the simultaneous computation of multiple multiscale basis functions. To overcome the discontinuities commonly observed in traditional activation functions, we further employ smoother activation function for the nonlinear transformations in different hidden layers. Numerical experiments demonstrate that the proposed method achieves the desired efficiency and accuracy (details are provided in \hyperref[sec:discussion]{Section \ref{sec:discussion}}).

The remainder of this paper is organized as follows. \hyperref[sec:preliminaries]{Section \ref{sec:preliminaries}} introduces the modeling background and the construction of multiscale basis functions. \hyperref[sec:method]{Section \ref{sec:method}} presents the architecture of the proposed deep learning framework. \hyperref[sec:result]{Section \ref{sec:result}} provides detailed numerical examples, including dataset construction, training strategy, and experimental results. \hyperref[sec:discussion]{Section \ref{sec:discussion}} discusses the findings in depth, and \hyperref[sec:conclusion]{Section \ref{sec:conclusion}} concludes our work.

\section{Preliminaries}
\label{sec:preliminaries}

\subsection{Model Problem and Finite Element Discretization}
Our problem begins with the following Darcy model for pressure field $p$ on a bounded Lipschitz computational domain $\Omega\in\mathbb{R}^2$:
\begin{equation}
	\begin{cases}
		-\text{div}(\kappa\nabla p)  &= f,~\text{in}~\Omega \\
		\kappa\nabla p \cdot \mathbf{n} &= 0,~\text{on}~\partial\Omega
	\end{cases}
	\label{eq:model_problem}
\end{equation}
where $\partial\Omega$ is a Lipschitz continuous boundary and $\mathbf{n}$ is the unit outward normal vector to $\partial\Omega$, $\kappa$ is the high-contrast permeability field with highly heterogeneous features, $f$ is given source term that satisfies the compatibility condition $\int_\Omega f~dx=0$, and there exists restriction $\int_\Omega p~dx=0$ to ensure the uniqueness of the solution. Induce the flux variable
\begin{equation}
	\mathbf{u}=-\kappa\nabla p \notag
\end{equation}
to the equation and apply the no-flux boundary condition, \eqref{eq:model_problem} can be transformed into the first-order form:
\begin{equation}
	\begin{cases}
		\kappa^{-1}\mathbf{u}+\nabla p &= 0,~\text{in}~\Omega~~~~(\text{Darcy's Law}) \\
		\text{div}(\mathbf{u}) &= f,~\text{in}~\Omega~~~~(\text{Mass Conservation}) \\
		\mathbf{u} \cdot \mathbf{n} &= g,~\text{on}~\partial\Omega~(\text{No-flux Boundary Condition})
	\end{cases}
\end{equation}
where $g$ is a given normal component of Darcy velocity on $\partial\Omega$.

We define $\mathcal{E}^H\coloneqq\bigcup_{i=1}^{N_e}E_i$ and $\mathcal{E}^h\coloneqq\bigcup_{i=1}^{M_e}$ as the set of all edges in the coarse and fine grids $\mathcal{T}^H$ and $\mathcal{T}^h$, where $N_e$ and $M_e$ are the number of coarse and fine grids, respectively.

Define
\begin{equation}
	L^2(\Omega)=\{ v:~v~\text{is defined in}~\Omega~\text{and square integrable, i.e.}~\int_\Omega v^2~dx~<~\infty \} \notag
\end{equation}
and use the Hilbert space $H(\text{div},\Omega)=\{ \mathbf{v}=(\mathbf{v}_1,~\mathbf{v_2})\in (L^2(\Omega))^2 \}$ and define
\begin{equation}
	V=H(\text{div},\Omega),~~W=L^2(\Omega) \notag
\end{equation}

We take the mixed finite element spaces on quaarilaterals:
\begin{equation}
	\begin{aligned}
		V_h &= \{ \mathbf{v}_h\in V: \mathbf{v}_h|_t = (b_tx_1+a_t, d_tx_2+c_t),~a_t,b_t,c_t,d_t\in\mathbb{R},~t\in\mathcal{T}^h \} \notag \\
		W_h &= \{ {w_h \in W: w_h~\text{is constant on each element in}~\mathcal{T}^h} \} \notag
	\end{aligned}
\end{equation}

Let $\{ \phi_i \}$ denote the set of multiscale basis functions for the coarse element, and the multiscale space for $p$ is defined as the span of all local basis functions:
\begin{equation}
	W_H=\text{span}\{\phi_i\} \notag
\end{equation}

Then the purpose of MGMsFEM is to find $(\mathbf{u}_H,p_H) \in (V_H, W_H)$ such that
\begin{equation}
	\begin{aligned}
		\int_\Omega \kappa^{-1}\mathbf{u}_H\cdot\mathbf{v}_H - \int_\Omega \text{div}(\mathbf{v}_H)p_H &= 0,~\quad\quad~~\forall~\mathbf{v}_H\in V_h^0 \\
		\int_\Omega \text{div}(\mathbf{u}_H)w_H &= \int_\Omega fw_H,~\forall~w_H \in W_H
	\end{aligned}
	\label{eq:mgmsfem}
\end{equation}
where $V_h^0=\{\mathbf{v}_h\in V_h : \mathbf{v}_h\cdot\mathbf{n}=0~\text{on}~\partial\Omega\}$. Let $\{\phi\}_{i=1}^n$ and $\{q_j\}_{j=1}^m$ denote the basis of $V_h$ and $W_h$, and assume that
\begin{equation}
	\mathbf{v}_h=\sum_{i=1}^n \mathbf{v}_i\phi_i,\quad p_h=\sum_{j=1}^m p_jq_j \notag
\end{equation}
then \eqref{eq:mgmsfem} can be written as the matrix representation:
\begin{equation}
	\begin{bmatrix}
		M & B^T \\
		B & 0
	\end{bmatrix}
	\begin{bmatrix}
		\mathbf{u} \\ p
	\end{bmatrix}
	=
	\begin{bmatrix}
		0 \\ F
	\end{bmatrix}
\end{equation}
where $M$ is a symmetric positive definite matrix with $M_{ij}=\int \kappa^{-1}\phi_i\phi_j$, $B$ is an approximation to the divergence operator $B_{ij}=-\int q_i\text{div}(\phi_j)$, and $F$ is a vector with $F_i=-\int fq_i$.

\subsection{Multiscale Basis Functions Construction}
After demonstrating the problem and total method of our numerical method, we will talk about the construction of multiscale basis functions (the core of our method) in this section.

Before computing the value of basis functions, it is necessary to construct the snapshot space, denoted by $W_{\text{snap}}$. In MGMsFEM, there are three ways to finish this step, and we will show them in detail here.

The first one is to take $W_{\text{snap}}$ as the fine-grid space for $W_h$, i.e.
\begin{equation}
	W_{\text{snap}} \coloneqq \{ \psi^{\text{snap}} \in W:~\text{piecewise constant on each coarse block} \} \notag
\end{equation}

The second is to solve local problems (including local spectral problems (LSPs) and local cell problems (LCPs)) with Dirichlet boundary condition:
\begin{equation}
	\begin{cases}
		\kappa^{-1}\mathbf{u}_j^{(i)}+\nabla p_j^{(i)} &= 0,~\text{in}~T_i \\
		\text{div}(\mathbf{u}_j^{(i)}) &= 0,~\text{in}~T_i \\
		p_j^{(i)}=\delta_j^{(i)} &= 
		\begin{cases}
			1,~\text{in}~e_i \\
			0,~\text{on other fine-edges on}~\partial T_i
		\end{cases},~j=1,\cdots,J_i
	\end{cases}
\end{equation}
where $J_i$ is the number of the element edges in the coarse block boundary. This problem can be solved numerically on the fine grid $T_i$ using lowest Raviart-Thomas element ($\text{RT}_0$ element) such that $p_j^{(i)}\in W_h$.

The last one is to solve local problems with Neumann boundary condition:
\begin{equation}
	\begin{cases}
		\kappa^{-1}\mathbf{u}_j^{(i)}+\nabla p_j^{(i)} &= 0,~~~~\text{in}~T_i \\
		\text{div}(\mathbf{u}_j^{(i)}) &= \alpha_j,~~\text{in}~T_i \\
		\frac{\partial p_j^{(i)}}{\partial \mathbf{n}_i} &= \delta_j^{(i)},~\text{on}~\partial T_i
	\end{cases}
\end{equation}
where $\mathbf{n}_i$ is the outward unit normal vector to $\partial T_i$, $\alpha_j$ satisfies satisfies the compatibility condition.
Hence, the snapshot space can be constructed by solving above problems and we have
\begin{equation}
	W_{\text{snap}} = \text{snap}\{ \psi_j^{\text{snap}},~j=1,\cdots,J_i,~\forall T_i \in \mathcal{T}^H \} \notag
\end{equation}

Next, we need to build the offline space. For each coarse element in $W_{\text{snap}}$, we reduce the spatial dimension by a LSP:
\begin{equation}
	\begin{aligned}
		a_i(p,w) &= \lambda s_i(p,w),~\forall~w\in W_{\text{snap}}^{(i)} \\
		a_i(p,w) &= \sum_{e\in\mathcal{E}_h^0} \kappa_e[p_h]_e[w_h]_e \\ 
		s_i(p,w) &= \int_{T_i} \kappa pw
	\end{aligned}
\end{equation}
where $[\cdot]$ is the jump, $e$ is an interior fine-scale edge in $T_i$. Then the discretization form of spectral problem can be written by
\begin{equation}
	A_iZ_k=\lambda_kS_iZ_k
\end{equation}
where $A_i$ and $S_i$ refer to the stiffness matrix and mass matrix, respectively. $\lambda_k$ is the eigenvalue and $Z_k$ is the corresponding eigenvector, arranged in ascending order. The first $l_i$ eigenvalues corresponding to the smallest eigenvectors are then selected to obtain the offline basis functions:
\begin{equation}
	\psi_k^{\text{off}}=\sum_{j=1}^{l_i} Z_{k,j}\psi_j^{\text{snap}}
\end{equation}
$Z_{k,j}$ represents the corresponding component of $Z_k$. The offline basis functions of all relevant element are combined to construct the global offline space:
\begin{equation}
	W_{\text{off}} \coloneqq \text{snap}\{ \psi_m^{\text{off}},~m=1,\cdots,M^{\text{off}},~M^{\text{off}}=\sum_{T_i\in\mathcal{T}^H}l_i \} \notag 
\end{equation}

It is worth noting that as the number of multiscale basis functions involved in the assembly of the stiffness matrix increases, the accuracy of the solution improves, but this also leads to higher computational cost. 

\section{Dual-Domain Deep Learning Framework}
\label{sec:method}
In this section, we will demonstrate the network structure of our method, and analyze its property on stability and convergence.

\subsection{Domain Transformation-based Feature Extraction}
\begin{figure}[h]
	\centering
	\includegraphics[width=1.0\linewidth]{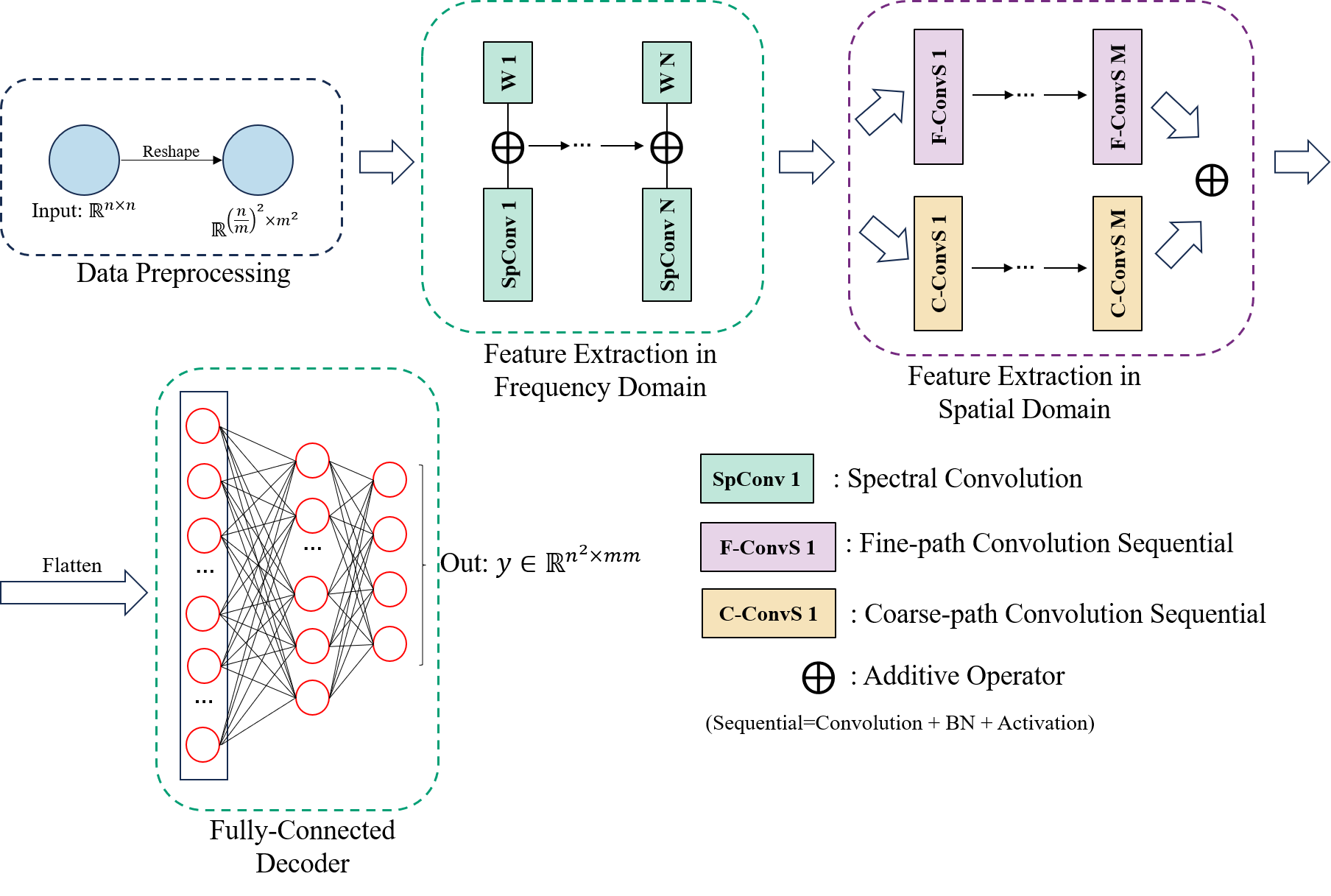}
	\caption{Network structure for our deep learning method. It can be divided into four parts: Data preprocessing, Feature Extraction in two domains, and Fully-Connected Decoder. Specially, the part ‘Feature Extraction in Frequency Domain’ also refers to FNO with TeLU function. The relevant abbreviations are explained in the lower right corner of the picture.}
	\label{fig:network}
\end{figure}

Our method origins from FNO, aiming at extracting frequency information inherent in the data. This part hires the Fourier integral operator $\mathcal{K}$ to transform the data from spatial domain to frequency domain, into a spectral representation.

Define the Fourier integral operator
\begin{equation}
	\mathcal{K}(\phi)v_t \coloneqq \mathcal{F}^{-1}(R_{\phi} \cdot (\mathcal{F}v_t)(x))
\end{equation}	
where $R_{\phi}$ represents the Fourier periodic function parameterized by $\phi \in \Theta_{\mathcal{K}}$. Here $\mathcal{K}$ plays a pivotal role in the domain transformation operation. We use the Fourier transform $\mathcal{F}$ to transfer the input data between different domains, which can be written as
\begin{equation}
	(\mathcal{F}K)(\xi)=\left\langle K,\psi \right\rangle_{L(D)}=\int_x v(x)\psi(x,\xi)\mu(x) \approx \sum_{x \in \mathcal{T}} v(x) 
\end{equation}
where $\psi(x,\xi)=\exp(2\pi i \left\langle x, \xi \right\rangle) \in L(D)$ is the Fourier basis function, $\xi$ refers to the frequency mode used for nonlinear transformations within the network, indicating the Fourier modes to be retained, $\mathcal{T}$ is a uniform grid that sampled from distribution $\mu$.

Operator $\mathcal{F}$ in this study refers to the two-dimensional fast Fourier transform (2D FFT). This serves as the first step, also the critical one in our network. In detail,
\begin{equation}
	\hat{x}(\xi) \coloneqq \int_{\Omega} x \cdot \exp(-2\pi i \left\langle x, \xi \right\rangle)
\end{equation}

In the frequency domain, the convolution operator can be converted into the elementwise multiplications, which is called spectral convolution. This process is a linear transformation in such domain. Through this method, we can obtain the frequency characteristics of our data and enhance the computation efficiency.

Let
\begin{equation}
	\hat{y}(\xi) \coloneqq \hat{\mathcal{K}}(\hat{\xi}) \cdot \hat{x}(\xi)
\end{equation}
where $\xi=(\xi_1,\xi_2)$, $\hat{\mathcal{K}}(\cdot)$ is a learnable kernel in the frequency domain. After spectral convolution, the modified data will be mapped back to the spatial domain via the inverse Fourier transform (iFFT):
\begin{equation}
	y(x) = \int_{\xi} \hat{y}\exp(2\pi i \left\langle x, \xi \right\rangle)~d\xi
\end{equation}

This process enables the model to recover spatial patterns after mode truncation, achieving the re-moval of redundant information while performing efficient feature extraction.

\subsection{Activation Function}
Classic FNO and convolutional neural networks (CNNs) use activation function that is not smooth enough (e.g. Rectified Linear Unit, ReLU). In this study we replace it with a smoother one, Tanh-exponential Linear Unit (TeLU), which is reported stable performance in many networks such as FNO and Deep Residual Net (ResNet) \cite{fernandez2024telu}. This can be written as
\begin{equation}
	\text{TeLU}(x)=x \cdot \tanh(e^x)
\end{equation}

This continuous function asymptotically approaches zero on the negative half-axis of x, while it exhibits an approximately linear behavior on the positive half-axis. This part can be seen in the part 'Feature Extraction in Frequency Domain' of \hyperref[fig:network]{Figure \ref{fig:network}}.

\subsection{Spatial Feature Extraction and Fusion}
After the feature extraction in the frequency domain, we attempt to obtain extra information from the spatial domain. Conventional neural network includes convolution layer, pooling layer, fully connected (linear) layer, and output layer. The number of parameters in a network generally in-creases with the complexity of its architecture, a trend that is particularly pronounced in the pres-ence of fully connected layers. To reduce complexity, one should minimize the excessive reliance on fully connected layers.

Based on this idea, we employ convolutional layers for feature extraction. Since multiscale basis functions are defined on the coarse grid, they can be expressed in terms of fine-grid information. To this end, we use convolutional kernels with sizes comparable to both coarse grids and fine grids to perform spatial feature extraction, which we denote as the coarse-grid path and the fine-grid path, respectively. Inspired by the concept of oversampling\cite{hou2009multiscale}, we consider feature extraction by reshaping the entire input field according to the coarse grid and concatenating the resulting matrices into a neural network input. This operation enables the net-work to simultaneously capture features across multiple grid scales. In other words, suppose the permeability field is given as $\kappa \in \mathbb{R}^{n \times n}$, and the coarse-grid block size is $m \times m$. Then the fine grid can be partitioned into $(\frac{n}{m})^2$ coarse blocks. After reshaping and rearranging, the resulting input matrix for the neural network has dimensions $\mathbb{R}^{(\frac{n}{m})^2 \times m^2}$

In both coarse and fine-grid path, the output data from FNO will be processed through a series of convolution sequential, which can be written as
\begin{equation}
	x^{(i+1)}=\sigma(\text{BN}(W^{(i)}\cdot x^{(i)}))
\end{equation}
Here, $W^{(i)}$ is the convolution filter, 'BN' refers to batch normalization and $\sigma$ is TeLU activation function. At the end of both paths, we leverage the additive operator to fuse the information:
\begin{equation}
	x_{\text{fused}} = x_{\text{coarse}} + x_{\text{fine}}
\end{equation}

This design enables us to capture spatial dependencies at different scales, effectively learning the multiscale information. After that, the data will be processed by multiple fully connected layers until the output fits the output size. To prevent the overfitting, we use ridge regression to give constraints:
\begin{equation}
	y=W_{\text{ridge}} \cdot x_{\text{fused}} + b_{\text{ridge}}
\end{equation}
with $L^2$ regularization loss
\begin{equation}
	\mathcal{L}=\frac{1}{N}\sum_{i=1}^{N} (y_i - \hat{y}_i)^2 + \lambda\| W_{\text{ridge}}\|_2^2
\end{equation}

The final output $y \in \mathbb{R}^{n^2 \times mm}$, where $mm$ is the used numbers of multiscale basis functions. For this, it can be seen in Figure \ref{fig:network}-'Feature Extraction in Spatial Domain' and 'Fully-Connected Decoder'.

\subsection{Analysis}
Let $\psi$ and $\hat{\psi}$ denote the multiscale basis functions from numerical and deep learning method, respectively. The LSP has been shown in the previous text. When the permeability field changes from $\kappa_1$ to $\kappa_2$, the pertubation of $a_i$ is
\begin{equation}
	\delta a_i(p,w) = \int_{T_i} (\kappa_1-\kappa_2) \nabla p \cdot \nabla w~dx
\end{equation}

By the perturbation theory of linear operator\cite{kato2013perturbation}, the change of eigenfunctions should satisfy
\begin{equation}
	\| \psi^{i,\text{off}}(\kappa_1) - \psi^{i,\text{off}}(\kappa_2)\|_{H^1(T_i)} \leq \frac{\| \kappa_1 - \kappa_2 \|_{L^{\infty}(T_i)}}{\kappa_{\text{min}}}\| \nabla\psi^{i,\text{off}}(\kappa_2)\|_{L^2(T_i)}
\end{equation}
where $\psi^{i,\text{off}}(\cdot)$ is offline basis function on $T_i$. Besides, according to the energy estimation of spectral problems
\begin{equation}
	\| \nabla\psi^{i,\text{off}}(\kappa) \|_{L^2} \leq \sqrt{\lambda^{(i)}_{\text{max}}} \| \psi^{i,\text{off}}(\kappa) \|_{L^2}
\end{equation}
and the normalization condition $\| \psi^{i,\text{off}}(\kappa)\|_{L^2}=1$, we have
\begin{equation}
	\| \psi(\kappa_1)-\psi(\kappa_2) \|^2_{H^1(T_i)} \leq L_{\text{local}} \|\kappa_1 - \kappa_2 \|_{L^{\infty}(T_i)},\quad L_{\text{Local}}=\frac{\sqrt{\lambda^{(i)}_{\text{max}}}}{\kappa_{\text{min}}}
\end{equation}
where $\kappa_{\text{min}}=\min \kappa$, $\kappa_{\text{min}} \leq \kappa \leq \kappa_{\text{max}}$. Then the difference of basis functions can be decomposed as
\begin{equation}
	\| \psi(\kappa_1) - \psi(\kappa_2) \|^2_{H^1(\Omega)} = \sum_{i=1}^{N_e} \| \psi^{i,\text{off}}(\kappa_1) - \psi^{i,\text{off}}(\kappa_2) \|^2_{H^1(T_i)}
	\label{eq:decompose}
\end{equation}

For each single term of the right side of \eqref{eq:decompose},
\begin{equation}
	\| \psi^{i,\text{off}}(\kappa_1) - \psi^{i,\text{off}}(\kappa_2) \|^2_{H^1(T_i)} \leq L^2_{\text{local}}\| \kappa_1 - \kappa_2 \|_{L^{\infty}}
\end{equation}

Sum up for all coarse blocks:
\begin{equation}
	\| \psi(\kappa_1) - \psi(\kappa_2) \|_{H^1(\Omega)} \leq L_{\text{global}}\sum_{i=1}^{N_e} \| \kappa_1 - \kappa_2 \|_{L^{\infty}(T_i)} \leq L_{\text{global}} \| \kappa_1 - \kappa_2 \|_{L^{\infty}(\Omega)}
\end{equation}
where $L_{\text{global}}=\max_i L_{\text{local}}^{(i)}$. Next, take the $H^1$-norm and use the triangle inequality, we have
\begin{equation}
	\| \hat{\psi}(\kappa_1) - \hat{\psi}(\kappa_2) \|_{H^1(\Omega)} \leq \underbrace{\| \hat{\psi}(\kappa_1) - \psi(\kappa_1) \|_{H^1(\Omega)}}_{\text{Term 1}} + \underbrace{\| \hat{\psi}(\kappa_2) - \psi(\kappa_2) \|_{H^1(\Omega)}}_{\text{Term 2}} + \underbrace{\| \psi(\kappa_1) - \psi(\kappa_2) \|_{H^1(\Omega)}}_{\text{Term 3}} \notag
\end{equation}

For the first two terms, it is easy to prove that $\forall~\varepsilon>0$, there exists that
\begin{equation}
	\| \hat{\kappa}_i - \psi(\kappa_i) \|_{H^1} \leq \varepsilon,\quad i=1,2 \notag
\end{equation}

This can be proved according to the Lipschitz continuity and the universal approximation theorem of neural operator\cite{cybenko1989approximation} and Sobolev embedding theorem. Substitute these terms into above inequality we can obtain
\begin{equation}
	\| \hat{\psi}(\kappa_1) - \hat{\psi}(\kappa_2) \|_{H^1} \leq 2\varepsilon + L_{\text{global}}\| \kappa_1 - \kappa_2 \|_{L^{\infty}}
	\label{eq:stability}
\end{equation}

\eqref{eq:stability} indicates that the difference of learned multiscale basis functions mostly comes from the difference of the input permeability fields, not from other instability.

Let the training set $\mathcal{D}$ be sampled from true distribution. Using above conclusion, the empirical error on $\mathcal{D}$ should converge uniformly to the polulation error as sampe size $N \rightarrow \infty$, i.e.,
\begin{equation}
	\lim_{N \rightarrow \infty} \sup_{\kappa \in \mathcal{D}} \| \hat{\psi} - \psi \|_{L^2} = 0 \notag
\end{equation}

For any $\varepsilon>0$, there exists a neural operator $\mathcal{N}$ and corresponding parameter set $\Theta^*$ such that
\begin{equation}
	\| \mathcal{N}(\kappa) - \hat{\mathcal{N}}(\kappa;\Theta^*) \|_{H^1} \leq \varepsilon_{\mathcal{N}}
\end{equation}

As $N \rightarrow \infty$, the empirical risk minimization ensures $\varepsilon_{\mathcal{N}} \rightarrow 0$. The stability of our network ensures the approximation is preserved under finite dataset size.

For the feature extraction in frequency domain, here we let $\mathcal{F}$ denote this part, it maps $\kappa$ to a low-dimensional spectral representation. From the spectral gap assumption, this part’s truncation error will decay exponentially, i.e.,
\begin{equation}
	\varepsilon_{\mathcal{F}} \leq C \cdot e^{-\gamma N_{\text{modes}}} \notag
\end{equation}
where $N_{\text{modes}}$ is the number of preserved Fourier modes. Then we have 
\begin{equation}
	\| \psi - \hat{\psi} \|_{L^2} \leq C_1 e^{-\gamma N_{\text{modes}}} + C_2 \varepsilon_{\mathcal{N}} \notag
\end{equation}
where $C_1$ and $C_2$ are constants. As $N_{\text{modes}},~N\rightarrow \infty$, both terms above will vanish:
\begin{equation}
	\lim\limits_{N_{\text{modes}},N\rightarrow\infty} \mathbb{E}_{\kappa\in\mathcal{D}}[\| \hat{\psi} - \psi \|_{L^2}] = 0
\end{equation}

\section{Numerical Examples}
\label{sec:result}
In this section, we will talk about the following topic: how we obtain the dataset, what strategy for model training, and the results we obtain.

\subsection{Data Generation and Dataset}
We leverage Karhunen-Loeve Expansion (KLE) as our permeability fields generator, which can generate stochastic fields with spatially correlated heterogeneity\cite{fukunaga1970application}. This approach offers a low-dimensional representation of the field via truncating high-order modes while maintaining dominant statistical features.

For a log-normal permeability field $\kappa(x;w)$, define the logarithmic transform
\begin{equation}
	\mathbb{Z}(x;w)=\log\kappa(x;w)
\end{equation}

This is modeled as a Gaussian random firld with mean $\mathbb{E}(x)$ and covariance $\mathbb{C}(x,y)$. KLE can decompose $\mathbb{Z}$ as
\begin{equation}
	\mathbb{Z}(x;w)=\mathbb{E}[x] + \sum_{i=1}^{\infty} \sqrt{\lambda_i}\phi_i(x)\xi_i(w)
\end{equation}
where $\lambda_i$ and $\phi_i$ are the eigenvalue and eigenfunction of the integral equation
\begin{equation}
	\int_{\Omega} \mathbb{C}\phi_i(y)dy = \lambda_i\phi_i(y)
\end{equation}
and $\xi_i(w)$ are the uncorrelated standard Gaussian random variables.

\begin{figure}[h]
	\centering
	\includegraphics[width=1.0\linewidth]{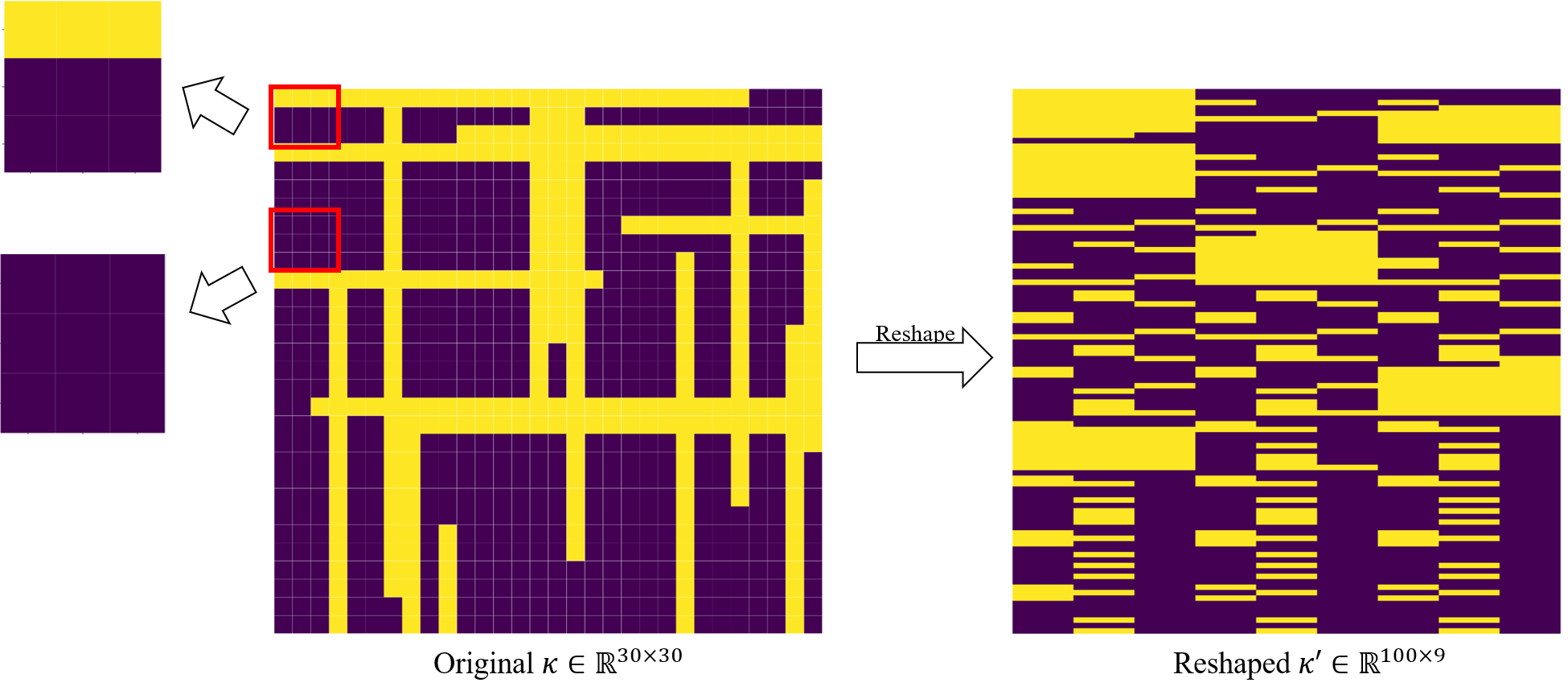}
	\caption{An example of our permeability field. Left: the field generated by KLE, with size $30 \times 30$. The yellow part represents the fissure, and the dark blue part represents the matrix. Right: reshaped permeability field. Each column refers to a coarse block. The coarse grid framed by the red line in the left figure is the coarse grid used for the subsequent display of multi-scale basis functions, one with cracks and one without.}
	\label{fig:data_sample}
\end{figure}

In this study, we choose the size of fine and coarse grid as $30 \times 30$ and $3 \times 3$, respectively. This means that each coarse block contains 9 fine grids, with total $10 \times 10$ coarse blocks. For this parameter setting, the system for single field have 100 coarse elements, where each element’s boundary consists of 12 fine-grid edges, a total of 1200 LCPs and 100 LSPs need to be solved, resulting in 1300 PDEs. Besides, as the first multiscale basis function is a piece-wise constant, which is not necessary to be trained by neural network, we choose 4 functions as our training targets (totally 5 functions, except the first one), i.e., the output $y \in \mathbb{R}^{900 \times 4}$. The data sample for this study can be seen in Figure \ref{fig:data_sample}.

The process of KLE and multiscale basis functions computation are complete via MATLAB. After the generation and duplication, we have totally 177800 samples, with 6537 duplicated samples. To construct the data loader for deep learning model training, we split them to three datasets: 102757 for training set, 34252 and 24354 for validation and test set, respectively. Meanwhile, each batch of data loader consists 64 randomly selected samples, which are transformed to tensor with size $64 \times 1 \times 100 \times 9$. This process is finished via Python--numpy and pytorch.

\subsection{Training and Optimization Strategy}
For parameter optimization, we employ the AdamW algorithm\cite{loshchilov2017decoupled}, a variant of Adam that decouples weight decay from the gradient-based updates, thereby improving generalization. The update rule can be expressed as
\begin{equation}
	\theta_{t+1} = \theta_t - \eta\frac{m_t}{\sqrt{v_t}+\epsilon} - \eta\lambda\theta_t
\end{equation}
where $\eta$ is the learning rate, $\lambda$ is the weight decay component, and $m_t,~n_t$ are the first- and second-moment estimates. Thanks to the automatic differential technique of Pytorch\cite{paszke2017automatic}, the gradients can be computed automatically, which applies the chain rule efficiently across all network layers.

To evaluate the performance of our deep learning method, we use mean squared error (MSE) and coefficient of determination ($\text{R}^2$) as metrics. Given the reference basis functions $y$ and predicted $\hat{y}$, the metrics are difined as
\begin{equation}
\begin{aligned}
	MSE_j&=\frac{1}{n}\sum_{i=1}^{n}(y_{ij}-\hat{y}_{ij})^2,~j=1,\cdots,mm \\
	\text{R}^2_j &= 1-\frac{\sum_{i=1}^n (y_{ij}-\hat{y}_{ij})}{\sum_{i=1}^n(y_{ij}-\bar{y}_j)^2}
\end{aligned} 
\end{equation}
where $mm$ is the number of multiscale basis functions and $\bar{y}_j$ is the mean of the j-th reference output. For model assessment, we report both the per-basis-function metrics and their average across all basis functions.

Note that the multiscale basis functions obtained by solving local problems are orthogonal to each other. In other words, given a coarse block, the corresponding multiscale basis function matrix $H \in \mathbb{R}^{9 \times 5}$ (here we include the first basis function) should satisfy
\begin{equation}
	\text{trace}(H^TH-I_{5 \times 5}) \rightarrow 0 \notag
\end{equation}

To verify whether the predicted basis function matrix of the model conforms to orthogonality, we define
\begin{equation}
	\text{Orth}=\frac{1}{N_e} \sum_{i=1}^{N_e} \| H_{i,\text{pred}}^T H_{i,\text{pred}} - I_{5 \times 5} \|_{L^2}^2
\end{equation}
where $H_{i,\text{pred}}$ denotes the predicted multiscale basis function matrix corresponding to the i-th coarse block in the entire permeability field. The smaller the Orth value, the better the orthogonality of predicted basis functions.

\subsection{Experiment Results}
In this section, we will illustrate our experiment environment and the results.

All neural network models were trained on an NVIDIA Tesla V100 GPU with 32 GB of memory. The implementation was based on PyTorch with CUDA acceleration. The learning rate was set to $1 \times 10^{-4}$, and the weight decay coefficient was $1 \times 10^{-3}$. Training was performed for a total of 300 epochs, and the best-performing model on the validation set was selected as the final model.

\begin{table}[h]
	\centering
	\caption{Evaluation metrics of our method: MSE and \( R^2 \) and Orth. This table includes the separate and average values of the four trained multiscale basis functions.}
	\begin{tabular}{lccccc}
		\toprule
		& Basis 2 & Basis 3 & Basis 4 & Basis 5 & Avg \\
		\midrule
		MSE       & 0.0019  & 0.0017  & 0.0009  & 0.0008  & 0.0013 \\
		\( \text{R}^2 \) & 0.9924  & 0.9931  & 0.9947  & 0.9891  & 0.9923 \\
		\midrule
		Orth      & \( 1.19 \times 10^{-3} \) &  &  &  &  \\
		\bottomrule
	\end{tabular}
	\label{tab:evaluation_metrics}
\end{table}

\begin{figure}[t]
	\centering
	\includegraphics[width=1.0\linewidth]{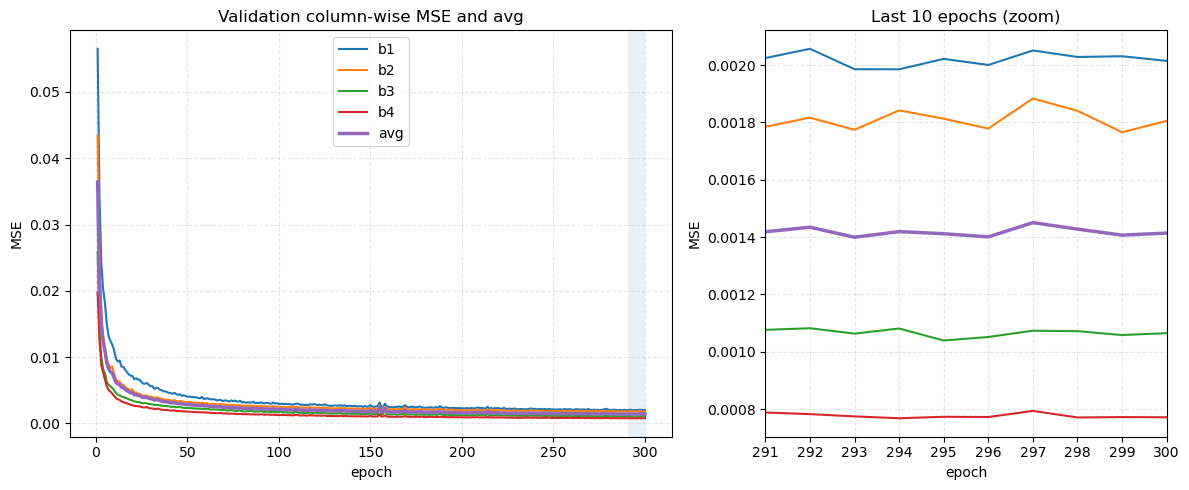}
	\label{fig:loss}
	\caption{ Learning curve of the training process. Left: the change of loss in total 300 epochs. Right: The loss curve of the last 10 epochs.}
\end{figure}

Table \ref{tab:evaluation_metrics} summarizes the evaluation metrics of the proposed method. The MSE for the trained basis functions remains on the order of $10^{-3}$, while the $\text{R}^2$ consistently exceeds 0.98, indicating high approximation accuracy. The value of Orth also indicates that the results of our method still follow the orthogonality. As shown in Figure \ref{fig:loss}, the training loss decreases rapidly during the initial epochs and stabilizes thereafter, with the loss curves of individual basis functions remaining nearly constant over the final 10 epochs, suggesting stable convergence.

\begin{figure}[h]
	\centering
	\includegraphics[width=1.0\linewidth]{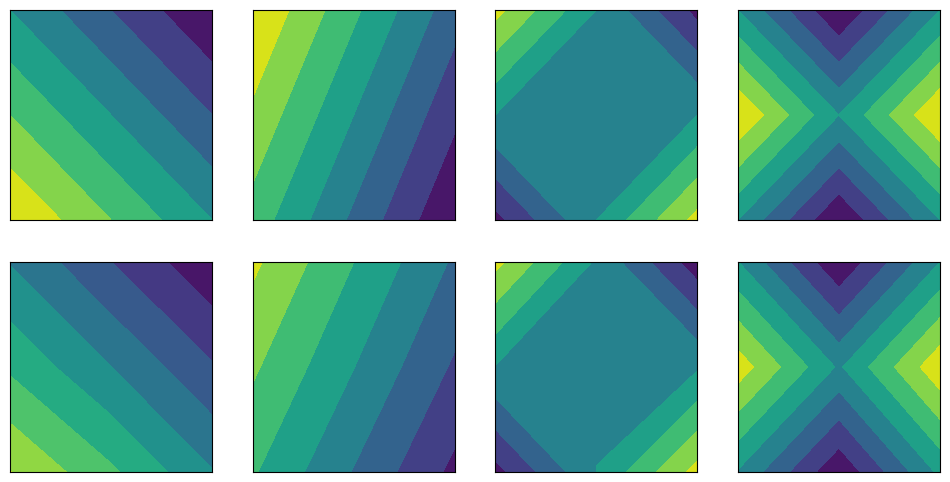}
	\label{fig:basis_matrix}
	\caption{The contour plots of our multiscale basis functions, which are for the coarse block with no fractures in Figure \ref{fig:data_sample}. Top: basis functions generated from MGMsFEM. Bottom: basis functions generated from our deep learning model.}
\end{figure}

\begin{figure}[h]
	\centering
	\includegraphics[width=1.0\linewidth]{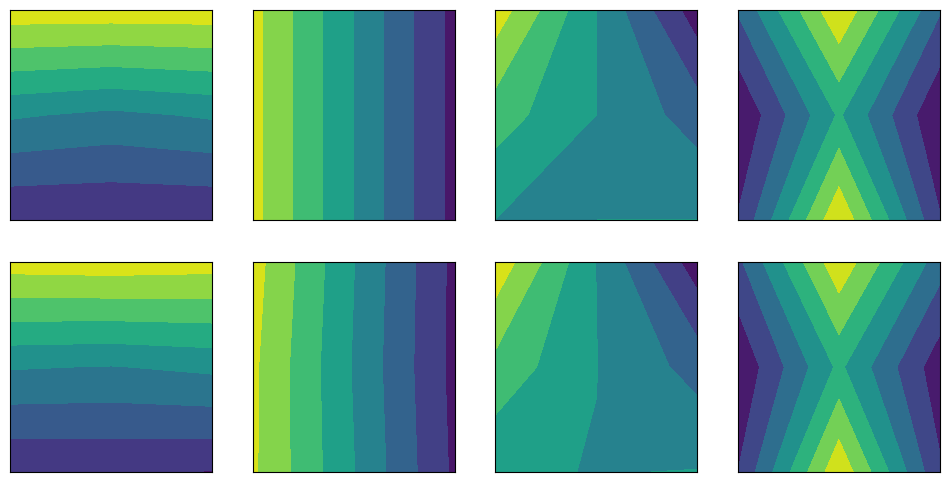}
	\label{fig:basis_fracture}
	\caption{The contour plots of our multiscale basis functions, which are for the coarse block with fractures in Figure \ref{fig:data_sample}. Top: basis functions generated from MGMsFEM. Bottom: basis functions generated from our deep learning model.}
\end{figure}

Figure \ref{fig:basis_matrix} and Figure \ref{fig:basis_fracture} present the contour plots of the multiscale basis functions. As shown in Figure \ref{fig:basis_matrix}, the basis functions on the coarse blocks of the matrix region exhibit a high level of agreement with the reference solutions, which is consistent with our expectations. In contrast, the contour plots in Figure \ref{fig:basis_fracture}, corresponding to coarse blocks containing fractures, reveal more noticeable discrepancies, particularly for Basis 4. Nevertheless, as a data-driven approach, our method is designed to approximate the underlying features based on large amounts of training data rather than to reproduce numerical solutions exactly. Such deviations are therefore expected. Importantly, when combined with the quantitative results reported in Table \ref{tab:evaluation_metrics}, these errors remain within a tol-erable range and do not compromise the overall effectiveness of the proposed framework.

\section{Discussion}
\label{sec:discussion}
In this study, we propose a dual-domain deep learning framework that can efficiently accelerate the computation of multiscale basis functions in heterogeneous porous media. The numerical experiments demonstrate that the MSE and Orth remain on the order of $10^{-3}$, while the $\text{R}^2$ consistently exceeds 0.98, indicating high approximation accuracy. The contour plots further confirm that the predicted basis functions in the matrix regions exhibit close agreement with the numerical solutions, whereas more noticeable discrepancies arise in coarse blocks containing fractures.

These discrepancies are reasonable and can be attributed to the inherent difficulty of representing highly heterogeneous features such as fractures. Unlike conventional numerical methods, which provide exact solutions to local boundary value problems, our approach is data-driven and relies on learning representative patterns from training data. Consequently, the method is designed to deliver fast and accurate approximations rather than exact solutions. Importantly, the quantitative metrics confirm that such deviations remain within a tolerable range and do not compromise the overall accuracy required for reservoir simulations.

Compared with traditional multiscale numerical methods such as MGMsFEM, the proposed framework significantly reduces computational cost while preserving accuracy. This improvement is particularly relevant for large-scale three-dimensional benchmark models such as SPE10, where the repeated construction of multiscale basis functions constitutes a major computational bottleneck. Moreover, in contrast to other deep learning approaches for PDEs—such as PINNs, and DeepONets - that primarily focus on approximating global solutions, our work targets the accelera-tion of basis function computation, filling a critical gap in the literature.

From the perspective of energy science, the ability to efficiently compute multiscale basis functions has practical implications for reservoir simulation, geological carbon sequestration, and groundwater resource management, all of which demand methods that balance accuracy with scalability. Future work will explore the integration of physics-informed constraints to further improve generalization, the extension of the framework to three-dimensional and multiphase flow problems, and the use of advanced architectures for enhanced interpretability. In addition, we well attempt to solve the final pressure field with a more efficient and accurate way based on current findings, and extend our framework to other cases.

\section{Conclusion}
\label{sec:conclusion}
In this work, we developed a dual-domain deep learning framework to accelerate the computation of multiscale basis functions in the MGMsFEM. By extracting features in both frequency and spatial domains and employing smoother activation functions, the method enables efficient and accurate construction of multiple basis functions simultaneously. Numerical experiments demonstrate that the proposed framework achieves a good trade-off between accuracy and efficiency, with MSE and Orth on the order of $10^{-3}$ and $\text{R}^2$ values exceeding 0.98, while significantly reducing computational cost. These findings indicate that integrating deep learning with multiscale numerical methods offers a promising direction for alleviating computational bottlenecks in reservoir simulation, and future work will extend the approach to three-dimensional and multiphase flow problems.

\section*{Author Contributions Statement}
\textbf{Peiqi Li}: Code, Methodology, Project Administration, Visualization, Writing-Original Draft; \textbf{Jie Chen}: Data Curation, Methodology, Supervision, Validation, Writing-Reviewing and Editing.

\section*{Conflict of Interest Statement}
The authors have no conflicts to disclose.

\section*{Code and Data Avaliability Statement}
The code and data that support the findings of this study are available from the corresponding author according to reasonable request.

\bibliographystyle{elsarticle-num-names}
\bibliography{ref}

\end{document}